\documentclass{amsart}
\usepackage{amsmath,amsthm,amssymb,amsfonts,amscd}
\usepackage{comment}

\newtheoremstyle{mytheorem}%
{}
{}
{\normalfont} 
{ } 
{\bfseries} 
{.} 
{ } 
{}

\theoremstyle{plain}
\newtheorem{Thm}{Theorem}[section]
\newtheorem{Prop}[Thm]{Proposition}

\newtheorem{Coro}[Thm]{Corollary}
\newtheorem{Lem}[Thm]{Lemma}
\newtheorem{Defn}[Thm]{Definition}

\theoremstyle{mytheorem}
\newtheorem{Remk}[Thm]{Remark}
\newtheorem{Ex}[Thm]{Example}

\newcommand{\G}{\mathcal{G}}
\newcommand{\gin}{\mathrm{gin}}
\newcommand{\iin}{\mathrm{in}}

\newcommand{\Tor}{\mathrm{Tor}}

\newcommand{\WLP}{weak Lefschetz property}

\newcommand{\SLP}{strong Lefschetz property}

\newcommand{\SSP}{strong Stanley property}
\newcommand{\SSE}{strong Stanley element}
\newcommand{\WSLP}{weak {\rm(}resp. strong{\rm)} Lefschetz property}
\newcommand{\WSLE}{weak {\rm(}resp. strong{\rm)} Lefschetz element}

\newcommand{\RA}{\rightarrow}

\numberwithin{equation}{section}


\begin{document}
\allowdisplaybreaks
\title[Generic Initial Ideal]
      {Generic Initial Ideals of Artinian ideals having Lefschetz Properties or the strong Stanley Property}

\author[Jea Man Ahn, Young Hyun Cho, and Jung Pil Park]
{Jea Man Ahn$^{\dag}$, Young Hyun Cho$^{\ddag}$, and Jung Pil
Park$^{\dag\dag}$}

\address{$\dag$ Computational Sciences,
                Korea Institute for Advanced Study,
                Seoul 130-722, South Korea }
\email{ajman@kias.re.kr} %

\address{$\ddag$ Department of Mathematical Sciences and Research Institute for Mathematics,
                Seoul National University,
                Seoul 151-747, South Korea }
\email{youngcho@math.snu.ac.kr} %

\address{$\dag\dag$ National Institute for Mathematical Sciences,
                 Daejeon 305-340, South Korea }
\email{jppark@nims.re.kr} %

\date{\today}

\thanks{The second author is partially supported by BK21, and the third author
is supported by National Institute for Mathematical Sciences.}

\begin{abstract}
 For a standard Artinian $k$-algebra $A=R/I$, we give equivalent conditions
for $A$ to have the weak (or strong) Lefschetz property or the
strong Stanley property in terms of the minimal system of generators
of the generic initial ideal $\mathrm{gin}(I)$ of $I$ under the
reverse lexicographic order. Using the equivalent condition for the
weak Lefschetz property, we show that some graded Betti numbers of
$\mathrm{gin}(I)$ are determined just by the the Hilbert function of
$I$ if $A$ has the weak Lefschetz property. Furthermore, for the
case that $A$ is a standard Artinian $k$-algebra of codimension 3,
we show that every graded Betti numbers of $\mathrm{gin}(I)$ are
determined by the graded Betti numbers of $I$ if $A$ has the weak
Lefschetz property. And if $A$ has the strong Lefschetz (resp.
Stanley) property, then we show that the minimal system of
generators of $\mathrm{gin}(I)$ is determined by the graded Betti
numbers (resp. by the Hilbert function) of $I$.
\end{abstract}

\maketitle

\section{\sc Introduction}
  Let $I$ be a homogeneous ideal of the polynomial ring $R=k[x_1,\ldots,x_n]$
over a field $k$. Throughout this paper, we use only a field $k$ of
characteristic 0. After Fr\"{o}berg \cite{Fr} introduced the
conjecture about the Hilbert series of generic algebras, it is one
of the remarkable problems under which condition $A=R/I$ has the
weak (or strong) Lefschetz property, which will be explicitly
defined in Definition $\ref{Def:WSLP}$. Recently, many authors
achieved some results on this problem. Watanabe \cite{Wa} showed
that "most" Artinian Gorenstein rings with fixed socle degree have
the \SLP(Cf. \cite{St}). In particular, it was proved by Watanabe
that any generic Artinian Gorenstein $k$-algebra has the \SLP. And
Harima, Migliore, Nagel and Watanabe \cite{HMNW} showed that $A$ has
the \SLP \ if $I$ is any complete intersection ideal of codimension
2, and has the \WLP \ if $I$ is any complete intersection ideal of
codimension 3. But it remains as an open problem whether any
complete intersection Artinian $k$-algebra $A$ has the \WLP \ or not
in higher codimension case.

  In this paper, we study Lefschetz properties and the strong Stanley property
in the view point of generic initial ideals. We begin Section 2 with
definitions of the weak (or strong) Lefschetz property and the
strong Stanley property of $A=R/I$. Then we introduce some results
of Wiebe who firstly investigated the Lefschetz properties in the
viewpoint of generic initial ideal. Wiebe \cite{Wie} proved that
$A=R/I$ has the \WSLP \ if and only if $R/\gin(I)$ has the \WSLP,
where $\gin(I)$ is the generic initial ideal of $I$ with respect to
the reverse lexicographic order (see Proposition
\ref{Prop:Equiv.Cond.I.and.GinI}). And he (or she) gave an
equivalent condition for $A=R/I$ to have the \WLP \ in terms of the
graded Betti numbers of $\gin(I)$ (see Proposition
\ref{Cond:BettiNumberForWLP}). We will give another equivalent
conditions and compute the graded Betti numbers of $\gin(I)$ more
precisely. To do so, we introduce the first reduction number and the
minimal system of generators of $\gin(I)$ in section 2.

  In section 3, we give equivalent conditions for $A$ to have the \WLP \ (Proposition
\ref{Prop:WLP}), the \SLP \ (Theorem \ref{Thm:SLP}), and the \SSP \
(Theorem \ref{Thm:SSP}), respectively, in terms of the minimal
system of generators of $\gin(I)$. And, as a result of Proposition
\ref{Prop:WLP}, we show that some graded Betti numbers of $\gin(I)$
can be computed just from the Hilbert function of $I$ (Corollary
\ref{Cor:BettiForWLP}).

  If we restrict ourselves to the case
$R=k[x_1,x_2,x_3]$, then we can obtain more information on the
minimal system of generators of $\gin(I)$ from the Hilbert function
or the graded Betti numbers of $R/I$. In section 4, under this
restriction for $R$, we show that if $R/I$ has the \WLP, then the
minimal generator $T$'s of $\gin(I)$ with $\max(T) \le 2$ are
uniquely determined by the graded Betti numbers of $R/I$
(Proposition \ref{Prop:UniquelyDetermineMinimalGeneratorWithWLP}).
Together with Corollary \ref{Cor:BettiForWLP}, this implies that
every graded Betti numbers of $\gin(I)$ can be computed from the
graded Betti numbers of $I$. And if $R/I$ has the \SLP, then we show
that the minimal system of generators of $\gin(I)$ is uniquely
determined by the graded Betti numbers of $R/I$ (Proposition
\ref{Prop:GinSLPUnique}). At last we show that if $R/I$ has the
\SSP, then the minimal system of generators of $\gin(I)$ is uniquely
determined by the Hilbert function of $R/I$ (Proposition
\ref{Prop:GinSSPUnique}).
\vskip 1cm%

\section{\sc Preliminaries}
  Let $R=k[x_1,x_2,\ldots,x_n]$ be the polynomial ring over a field
$k$. For a homogeneous ideal $I$ and a linear form $L$ in $R$, we
have the multiplication map
\begin{equation}\label{Mor:MultiMap}
  \times L^{i} : (R/I)_d \rightarrow (R/I)_{d+i},
\end{equation}
for each $i \ge 1$, $d \ge 0$.

\begin{Defn}
  A standard graded Artinian $k$-algebra
$A=R/I$ is said to have the \WSLP \ if there exists a linear form
$L$, called a \WSLE, such that the multiplication in
{\rm(}\ref{Mor:MultiMap}{\rm)} has maximal rank for each $d$ and
$i=1$ {\rm(}resp. for each $i \ge 1${\rm)}.
\end{Defn}

  Consider the exact sequence induced by $\times L^{i}$
\[
  0 \RA ((I:L^{i})/I)_{d} \RA (R/I)_{d} \xrightarrow[\qquad]{\times L^{i}}
         (R/I)_{d+i} \RA (R/(I+L^{i}))_{d+i} \RA
  0.
\]
If the multiplication $\times L^{i}$ has maximal rank, then either
$((I:L^{i})/I)_d$ or $(R/(I+L^{i}))_{d+i}$ must be 0 for each $d$.

There is more stronger property that a standard Artinian $k$-algebra
may have. That is the strong Stanley property named by Watanabe in
the paper \cite{Wa}, though it was originally called as the hard
Lefschetz property.

\begin{Defn}\label{Def:WSLP}
  A standard graded Artinian $k$-algebra
$A=R/I=\bigoplus_{i=0}^{t}A_i$, where $t = \max \{ i \ | \, A_i \neq
0 \}$, is said to have the \SSP \ if there exists a linear form $L
\in R$ such that the multiplication map $\times L^{t-2i}:A_i
\rightarrow A_{t-i}$ is bijective for each $i=0,1,\ldots,[t/2]$. In
this case the element $L$ is called a \SSE.
\end{Defn}

From the following exact sequence
\begin{equation}\label{ExSeq:SSP}
  0 \RA ((I:L^{t-2i})/I)_{i} \RA (R/I)_{i} \xrightarrow[\qquad]{\times L^{t-2i}}
         (R/I)_{t-i} \RA (R/(I+L^{t-2i}))_{t-i} \RA
  0,
\end{equation}
we can see that $A=R/I$ has the \SSP \ if and only if
$(I:L^{t-2i})_i = I_i$ and $R_{t-i} = (I+L^{t-2i})_{t-i}$ for all
$i=0,1,\ldots,[t/2]$.

\begin{Remk}\label{Remk:SymmetricSLP}
  A standard Artinian $k$-algebra $A=R/I$ has the \SSP \ if and only
if $A$ has the \SLP \ and the Hilbert function of $A$ is symmetric.
\end{Remk}

  In this paper, we will investigate those properties in the
viewpoint of the generic initial ideal of $I$ with respect to the
reverse lexicographic order. The notion of a generic initial ideal
originates with the following theorem of Galligo in \cite{Ga}.

\begin{Thm}\cite{Ga}
  For any multiplicative monomial order $\tau$ and any homogeneous
ideal $I \subset R$, there is a Zariski open subset $U \subset
GL_{n}(k)$ such that initial ideals $\rm{in}_{\tau}$$(gI)$ are
constant over all $g \in U$. We define the generic initial ideal of
$I$, denoted by $\rm {gin}_{\tau}$$(I)$, to be
\[
\gin_{\tau} (I):={\rm in}_{\tau} (gI)
\]
for $g \in U$.
\end{Thm}
  The most important property of generic initial ideals is that they
are {\it Borel-fixed}. In characteristic $0$ case, note that a
monomial ideal $I$ is Borel-fixed if and only if $I$ is {\it
strongly stable}: i.e., if $T$ is a monomial,
\[ x_i T \in I  \Rightarrow x_j T \in I, \qquad \forall
j\le i.
\]

  For a monomial ideal $I$, we denote by $\G(I)$ the
minimal system of generators of $I$, and by $\G(I)_d$ the set of
minimal generators of $I$ of degree $d$. For a monomial $T$ of $R$,
we set
  \[
    \max(T) = \max \{i \ | \, x_i \text{ divides } T \}.
  \]

 Eliahou-Kervaire Theorem gives us an easy way to compute the graded Betti
numbers of a stable monomial ideal $I$ of $R$.

\begin{Thm}[Eliahou-Kervaire] \cite{EK}
  Let $I$ be a stable monomial ideal of $R$. Then the graded Betti
number $\beta_{q,i}(I)=\dim_k{\rm Tor}^R_{q}(I,k)_i$ of $I$ is given
by
  \[
    \beta_{q,i}(I) = \sum_{T \in \G(I)_{i-q}} \binom{\max(T)-1}{q},
  \]
for each integers $q$ and $i$.
\end{Thm}

  Wiebe studied the Lefschetz properties in the view point of generic
initial ideals, and obtained the following results.

\begin{Lem}\cite{Wie}
  If $R/I$ is an Artinian ring such that $I$ is Borel-fixed,
  then the following conditions are equivalent:
  \begin{enumerate}
    \item $R/I$ has the \WSLP.
    \item $x_n$ is a \WSLE \ on $R/I$.
  \end{enumerate}
\end{Lem}

\begin{Prop}\cite{Wie}\label{Prop:Equiv.Cond.I.and.GinI}
  Let $R/I$ be a standard Artinian $k$-algebra, and
$\gin(I)$ the generic initial ideal of $I$ with respect to the
reverse lexicographic order. Then $R/I$ has the \WSLP \ if and only
if $R/\gin(I)$ has the \WSLP.
\end{Prop}

  Consider the following exact sequence
\begin{align}\label{ExSeq:LPGin}
\begin{split}
    0 \RA ((\gin(I):x_n^i)/\gin(I))_{d} \RA (R/\gin(I)&)_{d}
      \xrightarrow[\quad]{\times x_n^i} (R/\gin(I))_{d+i} \\
      & \RA (R/(\gin(I)+x_n^i))_{d+i} \RA 0.
\end{split}
\end{align}
  Then we can see that $R/I$ has the \WSLP \
if and only if either $(\gin(I):x_n^i)_d = \gin(I)_d$ or $R_{d+i} =
(\gin(I)+x_n^i)_{d+i}$ for each $d$ and $i=1$ (resp. for any $i \ge
1$).

  The following proposition is an equivalent condition for $R/I$ to have the \WLP \
in terms of the graded Betti numbers of $\gin(I)$.

\begin{Prop}\cite{Wie}\label{Cond:BettiNumberForWLP}
  Let $R/I$ be a standard Artinian $k$-algebra, and
$\gin(I)$ the generic initial ideal of $I$ with respect to the
reverse lexicographic order. If $d$ is the minimum of all $j \in
\mathbb{N}$ with $\beta_{n-1,n-1+j}(\gin(I)) > 0$, then the
following conditions are equivalent:
\begin{enumerate}
   \item $R/I$ has the \WLP.
   \item $\beta_{n-1,n-1+j}(\gin(I)) = \beta_{0,j}(\gin(I)) $ for
          all $j > d$.
   \item $\beta_{i,i+j}(\gin(I)) =
          \binom{n-1}{i}\beta_{0,j}(\gin(I)) $ for all $j > d$ and all $i$.
\end{enumerate}
\end{Prop}

  In Proposition \ref{Prop:WLP}, we will give another equivalent
condition for $R/I$ to have the \WLP. And we give a tool to compute
the graded Betti numbers of $\gin(I)$ from the Hilbert function of
$R/I$, if $R/I$ has the \WLP, in Corollary \ref{Cor:BettiForWLP}.
For these, we need the notion of reduction numbers.

\begin{Defn}
  If $I$ is a homogeneous ideal of $R$ with $\dim R/I = s$, then we
  define the $i$-th reduction number of $R/I$ to be $r_i(R/I) = \min\{t \ | \, \dim_k(R/(I+J_i))_{t+1} = 0 \}$ for $i \ge
  s$,  where $J_i$ is an ideal generated by $i$ general linear forms in $R$.
\end{Defn}

The following theorem implies that reduction numbers of $R/I$ and
$R/\gin(I)$ are the same, if we use the reverse lexicographic order.
\begin{Thm}\cite{HT}
  For a homogeneous ideal $I$ of $R$, if $\gin(I)$ is a generic
  initial ideal of $I$ with respect to the reverse
  lexicographic order, then
  \[
    r_{i}(R/I) = r_{i}(R/\gin(I))
               = \min \{t \ | \, x_{n-i}^{t+1} \in \gin(I) \}.
  \]
\end{Thm}

  In what follows, we will use only the reverse lexicographic order
as a multiplicative monomial order. Suppose that $I$ is a
homogeneous ideal of $R$ such that $R/I$ is a Cohen-Macaulay ring
with $\dim R/I = n-r$. In the paper \cite{CCP}, $\G(\gin(I))$ is
completely determined by the positive integer $f_1$ and functions
$f_i:\mathbb{Z}^{i}_{\ge 0} \longrightarrow \mathbb{Z}_{\ge 0} \cup
\{\infty\}$ defined as follows:
   \begin{align}\label{f's}
   \begin{split}
   f_1=&\min \{t \ | \, x^t_1 \in \gin(I) \} \text{ and, }\\
   f_i(\alpha_1, \ldots, \alpha_{i-1}) =& \min \{ t \ | \,
   x^{\alpha_1}_1 \cdots x^{\alpha_{i-1}}_{i-1} x^{t}_i \in
   \gin(I) \},
   \end{split}
  \end{align}
 for each  $ 2 \le i \le r$.

\begin{Lem}\cite{CCP} \label{Lem:f'sProperty}
  Let $f_1,\ldots,f_r$ be defined as in $(\ref{f's})$. For $2 \le i
\le r$, suppose that $\alpha_1,\ldots,\alpha_{i-1}$ are integers
such that $0 \le \alpha_1 < f_1$ and $0 \le \alpha_{j} <
f_{j}(\alpha_1, \ldots, \alpha_{j-1})$ for each $j\le i-1$. Then we
have
  \begin{enumerate}
  \item $0<f_i(\alpha_1, \ldots, \alpha_{i-1}) <\infty $,
  \item $x^{\alpha_1}_1 \cdots x^{\alpha_{i-1}}_{i-1}
         x^{f_i(\alpha_1, \ldots, \alpha_{i-1})}_i \in
         \G(\gin(I))$, and
  \item if $1 \le j \le i-1$ and $\alpha_j \ge 1$, then
   \[f_i(\alpha_1, \ldots, \alpha_j, \ldots, \alpha_{i-1}) \le
         f_i(\alpha_1, \ldots, \alpha_j-1, \ldots, \alpha_{i-1}) - 1.
   \]
  \end{enumerate}
\end{Lem}

For $ 1 \le i \le r-1$, let
\begin{equation} \label{Set:IndexMinGen}
\begin{split}
  J_i = \left\{(\alpha_1,\ldots,\alpha_{i})
              \left| \begin{array}{l}
                       0 \le \alpha_1 < f_1, \ \text{and for each} \ 2 \le j \le i, \\
                       0 \le \alpha_{j} < f_{j}(\alpha_1, \ldots, \alpha_{j-1})
                     \end{array}
              \right.
        \right\},
\end{split}
\end{equation}
and let
\begin{equation} \label{Set:MinGen}
\begin{split}
   \G=\{x_1^{f_1}\} \cup
       \left\{
              x^{\alpha_1}_1 \cdots x^{\alpha_{i-1}}_{i-1}
              x^{f_i(\alpha_1, \ldots, \alpha_{i-1})}_i
              \left| \begin{array}{l}
                       (\alpha_1,\ldots,\alpha_{i-1}) \in J_{i-1}, \\
                       \text{for} \ 2 \le i \le r
                     \end{array}
              \right.
            \right\}.
\end{split}
\end{equation}

\begin{Prop}\cite{CCP}\label{Prop:MinimalGen}
   Let $f_1,\ldots,f_r$ be defined as in $(\ref{f's})$. Then
the minimal system of generators of $\gin(I)$ is $\G$ in
$(\ref{Set:MinGen})$.
\end{Prop}

\begin{Remk}\label{Cond:ElementInJ_i}
  \begin{enumerate}
  \item Let $1 \le i \le r-1$
        and $(\alpha_1,\ldots,\alpha_i) \in \mathbb{Z}_{\ge 0}^i$.
        Then $(\alpha_1,\ldots,\alpha_i)$ belongs to $J_i$ if and only if
        $x_1^{\alpha_1} \cdots x_i^{\alpha_i}
        x_{i+1}^{f_{i+1}(\alpha_1,\ldots,\alpha_i)}$ is an element of
        $\G(\gin(I))$ and $0 < f_{i+1}(\alpha_1,\ldots,\alpha_i) < \infty$
        by Lemma \ref{Lem:f'sProperty} (1).
  \item If $(\alpha_1,\ldots,\alpha_i) \in J_i$, then the element
$(0,\ldots,0,|\alpha|) \in \mathbb{Z}_{\ge 0}^{i}$ belongs to $
J_i$, where $|\alpha| := \sum_{j=1}^{i} \alpha_j$. In particular,
$|\alpha| \le f_i(0,\ldots,0) - 1$.
\begin{proof}
By the definition of $J_{i}$, it is enough to show that $|\alpha| <
f_{i}(0,\ldots,0)$. If $|\alpha|=0$, then the assertion follows,
hence we may assume that $|\alpha| \ge 1$. Suppose to the contrary
that $|\alpha| \ge f_{i}(0,\ldots,0)$. Then $x_{i}^{|\alpha|} \in
\gin(I)$. Since $\gin(I)$ is strongly stable, $x_1^{\alpha_1}\cdots
x_{i}^{\alpha_{i}} \in \gin(I)$. This contradicts to (1).
\end{proof}
  \item Let $(\alpha_1,\ldots,\alpha_s,\ldots,\alpha_t,\ldots,\alpha_i) \in J_i$,
        $1 \le i \le r-1$. If $\alpha_s \ge 1$, then $f_{i+1}(\alpha_1,\ldots,
        \alpha_s,\ldots,\alpha_t,\ldots,\alpha_i) \le f_{i+1}(\alpha_1,\ldots,
        \alpha_s - 1, \ldots, \alpha_t + 1, \ldots, \alpha_i)
        $,
        since $\gin(I)$ is strongly stable.
  \end{enumerate}
\end{Remk}

\begin{Remk}\label{Remk:JDetermineGin}
  Note that the set $J_{r-1}$ determines the minimal generator $T$'s
of $\gin(I)$ satisfying $\max(T) \le r-1$. Indeed, if $i \le r-1$
and $x_1^{\alpha_1}\cdots x_{i-1}^{\alpha_{i-1}} x_i^{\alpha_i} \in
\G(\gin(I))$ with $\alpha_i \ge 1$, then
$(\alpha_1,\ldots,\alpha_{i-1},\alpha_i-1,0,\ldots,0) \in
\mathbb{Z}_{\ge 0}^{r-1}$ belongs to $J_{r-1}$. Hence if we set
$g_i(\alpha_1,\ldots,\alpha_{i-1}) := \max \{ \beta \ | \,
(\alpha_1,\ldots,\alpha_{i-1},\beta,0,\ldots,0) \in J_{r-1} \} + 1$,
then $g_i(\alpha_1,\ldots,\alpha_{i-1}) = \alpha_i$. This shows that
the assertion follows.
\end{Remk}
  The following example will help to understand above remarks.
\begin{Ex}
  Suppose that an homogeneous ideal $I$ of $R=k[x,y,z,w]$ has the following
the generic initial ideal:
\[
\gin(I)=(x^2, xy, y^3, y^2z, xz^2, yz^2, z^3, xzw^2 , y^2w^3, yzw^3,
z^2w^3, xw^5, yw^5, zw^5, w^7).
\]
Since $f_1=2$, we have $J_1=\{ \, 0,1\}$.  And $f_2(0)=3$,
$f_2(1)=1$ imply that $J_2=\{(0,0),(0,1),(0,2),(1,0)\}$. And note
that
\[ J_3=\{(0,0,0),(0,0,1),(0,0,2),(0,1,0),(0,1,1),(0,2,0),(1,0,0),(1,0,1)\}, \]
since we have $f_3(0,0)=3$, $f_3(0,1)=2$, $f_3(0,2)=1$,
$f_3(1,0)=2$.
\end{Ex}

\vskip 1cm %

\section{\sc Equivalent Conditions}
  Throughout this section we assume that $A=\bigoplus_{i=0}^{t} A_i=R/I$ is a
standard Artinian algebra over a field $k$, where $I$ is a
homogeneous ideal of $R=k[x_1,\ldots,x_n]$ and $t=\max\{i \ | \, A_i
\neq 0 \}$.

  In Proposition \ref{Cond:BettiNumberForWLP}, an equivalent
condition for $A$ to have the \WLP \ was given in terms of the
graded Betti numbers of $\gin(I)$. In this section, we give
equivalent conditions for $A$ to have the \WLP \ (Proposition
\ref{Prop:WLP}), the \SLP \ (Theorem \ref{Thm:SLP}), and the \SSP \
(Theorem \ref{Thm:SSP}), respectively, in terms of the minimal
system of generators of $\gin(I)$. Furthermore, as a result of
Proposition \ref{Prop:WLP}, we show that graded Betti numbers of
$\gin(I)$ can be computed just from the Hilbert function of $I$ (See
Corollary \ref{Cor:BettiForWLP}).

  To see the relation between the Hilbert function of $I$ and
the minimal generator $T$'s of $\gin(I)$ with $\max(T) = n$,
consider the following lemma and proposition. Although these are
introduced in the paper \cite{AS}, we give full proofs of them
because we modified a little.

\begin{Lem}\label{Lem:Cond_Belong_In_G(ginI)}
  If $T$ is a nonzero monomial in
  $((\gin(I):x_n^i)/(\gin(I):x_n^{i-1}))_{d}$ for $d \ge 0$ and $i \ge 1$, then $Tx_n^i \in
  \G(\gin(I))_{d+i}$.
\end{Lem}
\begin{proof}
  Suppose to the contrary that $Tx_n^i$ is not a minimal generator
of $\gin(I)$. Then $Tx_n^i = uM$ for some monomials $M \in \gin(I)$,
$u \in R$ with $\deg(u) \ge 1$. Since $T \notin
(\gin(I):x_n^{i-1})$, $x_n$ cannot divide $u$, i.e. $\max(u) < n$,
and so $x_n$ divides $M$. But since $\gin(I)$ is strongly stable,
this implies  that
\[
    Tx_n^{i-1} = u \frac{M}{x_n} \in \gin(I),
\]
which is a contradiction.
\end{proof}

\begin{Prop}\label{Prop:r1A}
    Suppose that $d > r_{1}(A)$. If $T_1,\ldots,T_t$ are
monomials which form a $k$-basis of $((\gin(I):x_n)/\gin(I))_{d-1}$,
then
  \[
    \{ T_1x_n, \ldots, T_tx_n \} =
      \{ T \in \G(\gin(I)) \ | \, \max(T) = n, \deg(T) = d \}.
  \]
In particular,
   \begin{align}\label{Eq:d>r_1(A)}
   \begin{split}
       H(A,d-1) - H(A,d)
                =& \ \dim_k ((\gin(I):x_n)/\gin(I))_{d-1} \\
                =& \ \sharp \{ T \in \G(\gin(I)) \ | \, \max(T) = n, \deg(T) = d
\},
   \end{split}
   \end{align}
where $H(A,d) := \dim_k A_d$ is the Hilbert function of $A$ at
degree $d$. Thus we have $H(A,d-1) \ge H(A,d)$ for any $d >
r_{1}(A)$.
\end{Prop}
\begin{proof}
  Note that if $T \in \G(\gin(I))$ is a monomial with $\deg(T) = d$ and
$\max(T) = n$, then $T/x_n$ is a nonzero monomial in
$((\gin(I):x_n)/\gin(I))_{d-1}$. Hence, it is enough to show that if
$T_1,\ldots,T_t$ are monomials which form a $k$-basis of
$((\gin(I):x_n)/\gin(I))_{d-1}$, then $T_ix_n \in \G(\gin(I))$ for
all $i$. But this follows from Lemma
\ref{Lem:Cond_Belong_In_G(ginI)}.

  For the second assertion, note that
$d > r_1(A)$ implies $x_{n-1}^{d} \in \gin(I)$ by the definition of
the reduction number $r_{1}(A)$. Since $\gin(I)$ is strongly stable,
this implies that $(R/(\gin(I)+ x_n))_{d} = 0$. From the exact
sequence (\ref{ExSeq:LPGin}), we have
\begin{align*}
 \dim_k ((\gin(I):x_n)/\gin(I))_{d-1}
           =& \ H(R/\gin(I),d-1) - H(R/\gin(I),d) \\
           =& \ H(A,d-1) - H(A,d).
\end{align*}
The second equality in (\ref{Eq:d>r_1(A)}) follows from the first
assertion.
\end{proof}

\begin{Coro}\label{Cor:r_1WithWLP}
 Suppose that $A=R/I$ has the \WLP. If $H(A,d-1) \ge H(A,d)$,
 then $d > r_{1}(A)$. In particular,
       \[
          r_{1}(A) = \min \{ d \ | \, H(A,d-1) \ge H(A,d) \} - 1.
       \]
\end{Coro}
\begin{proof}
  If $d \le r_1(A)$, then $x_{n-1}^d \notin \gin(I)$. So $R_{d} \neq (\gin(I) + x_n)_{d}$.
Since $R/\gin(I)$ has the \WLP, we have $(\gin(I) : x_n)_{d-1} =
(\gin(I))_{d-1}$. By the exact sequence (\ref{ExSeq:LPGin}), this
implies that
\[
   H(A,d-1) = H(A,d) - \dim_k(R/(\gin(I) + x_n))_{d} < H(A,d),
\]
which is a contradiction. The last assertion follows from
Proposition \ref{Prop:r1A} and the first assertion.
\end{proof}

  The following proposition gives us a criterion to check that $A=R/I$
has the \WLP, if we know the minimal system of generators or the
graded Betti numbers of $\gin(I)$ and the first reduction number
$r_1(A)$. Remind $J_{n-1}$ defined in (\ref{Set:IndexMinGen}).

\begin{Prop}\label{Prop:WLP}
  The following statements are equivalent:
  \begin{enumerate}
    \item $A=R/I$ has the \WLP.
    \item If $T$ is a minimal generator of $\gin(I)$ of
          degree $d$ with $\max(T) = n$, then
          $d \ge r_{1}(R/I) + 1$, i.e.
          \[
             |\alpha| +  f_n(\alpha_1,\ldots,\alpha_{n-1}) \ge r_1(A) + 1,
          \]
          for any $(\alpha_1,\ldots,\alpha_{n-1}) \in J_{n-1}$,
          where $|\alpha| = \sum_{i=1}^{n-1} \alpha_i$.
    \item $\beta_{n-1,n-1+j}(\gin(I)) = 0$ for all $j \le r_1(A)$.
  \end{enumerate}
\end{Prop}
\begin{proof}
  (1) $\Rightarrow$ (2):
   If $T$ is a minimal generator of $\gin(I)$ of degree $d$ with
$\max(T) = n$, then  $T/x_n$ is a nonzero element in
$((\gin(I):x_n)/(\gin(I))_{d-1}$. Since $A$ has the \WLP,
$(R/(\gin(I) + x_n)_{d} = 0$. Hence $x_{n-1}^d
  \in \gin(I)$, this implies that $d \ge r_1(A) + 1$.

  (2) $\Rightarrow$ (1): Suppose that $((\gin(I):x_n)/\gin(I))_{d-1} \neq
0$. Let $T$ be a monomial in $(\gin(I):x_n)_{d-1}$ such that $T
\notin \gin(I)$. Then we have $Tx_n \in \G(\gin(I))_{d}$, by Lemma
\ref{Lem:Cond_Belong_In_G(ginI)}. By assumption, we have $d \ge
r_1(A) + 1$. This implies that $x_{n-1}^{d} \in \gin(I)$, so
$(R/(\gin(I)+x_n))_{d} = 0$.

  (2) $\Leftrightarrow$ (3): It is clear from Eliahou-Kervaire Theorem.
\end{proof}

\begin{Coro}\label{Cor:BettiForWLP}
  If $A=R/I$ has the \WLP, then
  \begin{enumerate}
  \item $\beta_{n-1,n-1+d}(\gin(I))=
           \begin{cases}
                0&                   \ \text{ if } d \le r_1(A), \\
                H(A,d-1)-H(A,d)& \ \text{ if } d > r_1(A).
           \end{cases}
        $
  \item For any $0 \le i \le n-1$ and $d \ge r_1(A) + 2$,
        \[
          \beta_{i,i+d}(\gin(I)) =
          \binom{n-1}{i}(H(A,d-1)-H(A,d)).
        \]
  \end{enumerate}
\end{Coro}
\begin{proof}
  The assertions follow from Proposition \ref{Prop:r1A}, Proposition
  \ref{Prop:WLP}, and Eliahou-Kervaire Theorem.
\end{proof}

\begin{Ex}
  Consider two strongly stable ideals of $R=k[x_1,x_2,x_3]$ defined by
  \begin{align*}
   I &= (x_1^2, x_1x_2, \underline{x_1x_3, x_2^3}, x_2x_3^2, x_3^4), \\
   J &= (x_1^2, x_1x_2, \underline{x_2^2, x_1x_3^2}, x_2x_3^2, x_3^4).
  \end{align*}
  Since $\deg(x_1x_3) = 2 < 3 = r_1(R/I) + 1$, $R/I$ does not have the \WLP, but $R/J$ has.
\end{Ex}

  The following theorem gives an equivalent condition for $A=R/I$ to have
the \SLP \ in terms of the minimal system of generators of
$\gin(I)$.

\begin{Thm}\label{Thm:SLP}
$A=R/I$ has the \SLP \ if and only if
 for any $\alpha=(\alpha_1,\ldots,\alpha_{n-1}) \in J_{n-1}$,
 \begin{enumerate}
 \item $|\alpha| + f_n(\alpha_1,\ldots,\alpha_{n-1})
       \ge r_1(A) + 1$, and
 \item $f_n(\alpha_1,\ldots,\alpha_{n-1}) \ge
        f_n(0,\ldots,0,|\alpha| + 1) + 1 $,
 \end{enumerate}
 where $|\alpha| = \sum_{j=1}^{n-1} \alpha_j$.
\end{Thm}
\begin{proof}
Suppose that $A$ has the \SLP, and $(\alpha_1,\ldots,\alpha_{n-1})
\in J_{n-1}$. The first assertion follows from Proposition
\ref{Prop:WLP}. For the second assertion, note that
$(0,\ldots,0,|\alpha|) \in J_{n-1}$ by Remark
\ref{Cond:ElementInJ_i} (2). Hence if $(0,\ldots,0,|\alpha|+1)
\notin J_{n-1}$, then $x_{n-1}^{|\alpha|+1} \in
 \gin(I)$ by the definition of $J_{n-1}$. So we have $f_n(0,\ldots,0,|\alpha|+1) =
0$, and the second assertion follows from Lemma
\ref{Lem:f'sProperty} (1). Thus we can assume that
$(0,\ldots,0,|\alpha|+1) \in J_{n-1}$. For simplicity, let $\beta =
f_n(\alpha_1, \ldots, \alpha_{n-1})$. Then $x_1^{\alpha_1} \cdots
x_{n-1}^{\alpha_{n-1}} x_n^{\beta} \in \G(\gin(I))$. This shows that
$x_1^{\alpha_1} \cdots x_{n-1}^{\alpha_{n-1}}$ is a nonzero monomial
in $((\gin(I):x_n^{\beta})/\gin(I))_{|\alpha|}$. From the \SLP \ of
$A$, we have $R_{|\alpha|+\beta} =
(\gin(I)+x_n^{\beta})_{|\alpha|+\beta}$, and hence
$x_{n-1}^{|\alpha|+1} x_n^{\beta-1} \in \gin(I)$. By the definition
of $f_n$, this implies that $f_n(0,\ldots,0,|\alpha|+1) \le \beta -
1$.

 For the converse, it is enough to show that if $((\gin(I) :
x_n^i)/\gin(I))_{d} \neq 0$, then $R_{d+i} = (\gin(I) +
x_n^i)_{d+i}$ for all $d$ and $i$. We will show this by induction on
$i$. By Proposition \ref{Prop:WLP}, $A$ has the \WLP, i.e.
$(R/(\gin(I) + x_n^i))_{d+i} = 0$ for all $d$ and $i=1$. Suppose
that $i > 1$ and that $T=x_1^{\alpha_1} \cdots
x_{n-1}^{\alpha_{n-1}} x_n^{\alpha_n}$ is a nonzero monomial in
$((\gin(I) : x_n^i)/\gin(I))_{d}$. If $T$ is a nonzero monomial in
$((\gin(I) : x_n^{i-1})/\gin(I))_{d}$, then $R_{d+i-1} = (\gin(I) +
x_n^{i-1})_{d+i-1}$, by induction hypothesis. So $x_{n-1}^{d+1}
x_n^{i-2} \in \gin(I)$, and hence $x_{n-1}^{d+1} x_n^{i-1} \in
\gin(I)$. Since $\gin(I)$ is strongly stable, this shows that
$R_{d+i}=(\gin(I) + x_n^i)_{d+i}$. Thus we may assume that $T$ is a
nonzero monomial in $((\gin(I) : x_n^i)/(\gin(I):x_n^{i-1}))_{d}$.
Then $Tx_n^i = x_1^{\alpha_1} \cdots x_{n-1}^{\alpha_{n-1}}
x_n^{\alpha_n + i} \in \G(\gin(I))$ by Lemma
\ref{Lem:Cond_Belong_In_G(ginI)}. This implies that
$(\alpha_1,\ldots,\alpha_{n-1}) \in J_{n-1}$ by Remark
\ref{Cond:ElementInJ_i}(1) and
\[
\alpha_n + i = f_n(\alpha_1,\ldots,\alpha_{n-1}) \ge
f_n(0,\ldots,0,d - \alpha_n + 1) + 1,
\]
since $d=\deg(T)=\sum_{j=1}^{n} \alpha_j$ and $f_n$ satisfies the
condition (2). By the definition of $f_n$, we have $x_{n-1}^{d -
\alpha_n + 1} x_n^{\alpha_n + i - 1} \in \gin(I)$. Since $\gin(I)$
is strongly stable, this shows that $x_{n-1}^{d + 1} x_n^{i-1} \in
\gin(I)$, and hence the assertion follows.
\end{proof}

\begin{Ex}\label{Ex:NonUniqueSLP}
Consider the following stable monomial ideal of $R=k[x,y,z,w]$
\[
I=(x^2, xy, y^3, y^2z, xz^2, yz^2, z^3, xzw^2, y^2w^3, yzw^3,
z^2w^3, xw^5, yw^5, zw^5, w^7).
\]
Every generator $T$ which is divisible by $w$ has degree $ \ge 3$,
and $r_1(R/I) = 2$. Furthermore, $I$ satisfies the second condition
in Theorem \ref{Thm:SLP}. Hence $R/I$ has the \SLP.
\end{Ex}

  Finally, we will give an equivalent condition for $A=R/I$ to have the
\SSP \ in terms of the minimal system of generators of $\gin(I)$. To
do so, note that $A=R/I$ has the \SSP \ if and only if $R/\gin(I)$
has the \SSP \ with a \SSE \ $x_n$, as shown in the cases of \WLP \
and \SLP \ (cf. Proposition \ref{Prop:Equiv.Cond.I.and.GinI}).
Consider the following exact sequence
\begin{align*}
    0 \RA ((\gin(I):x_n^{t-2i})/\gin(I))_{i} \RA (R/\gin(I)&)_{i}
      \xrightarrow[\quad]{\times x_n^{t-2i}} (R/\gin(I))_{t-i} \\
      & \RA (R/(\gin(I)+ x_n^{t-2i}))_{t-i} \RA 0.
\end{align*}
 Then we can see that $A=R/I$ has the \SSP \ if and only if
$(\gin(I):x_n^{t-2i})_i = \gin(I)_i$ and $R_{t-i} = (\gin(I)+
x_n^{t-2i})_{t-i}$, or equivalently $x_{n-1}^{i+1} x_n^{t-2i-1} \in
\gin(I)_{t-i}$ for any $i=0,1,\ldots,[t/2]$, since $\gin(I)$ is
strongly stable.

\begin{Lem}\label{Lem:SizeOfAlpha}
  Suppose that $A=R/I=\bigoplus_{i=0}^{t} A_i$ has the \SSP. If
$(\alpha_1,\ldots,\alpha_{n-1}) \in J_{n-1}$, then $|\alpha| :=
\sum_{i=0}^{n-1} \alpha_i \le [t/2]$.
\end{Lem}
\begin{proof}
Since the Hilbert function of $A$ is symmetric, $H(A,[t/2]) \ge
H(A,[t/2]+1)$. By Corollary \ref{Cor:r_1WithWLP}, this implies that
$r_1(A) \le [t/2]$. Since $(\alpha_1,\ldots,\alpha_{n-1}) \in
J_{n-1}$, we have%
\[
   |\alpha| \le f_{n-1}(0,\ldots,0) - 1 = r_1(A) \le [t/2],
\]
by Remark \ref{Cond:ElementInJ_i} (2).
\end{proof}

\begin{Thm}\label{Thm:SSP}
A standard Artinian $k$-algebra $A=R/I=\bigoplus_{i=0}^{t}A_i$ has
the \SSP \ if and only if for any $(\alpha_1,\ldots,\alpha_{n-1})
\in J_{n-1}$,
\begin{equation}\label{Eq:SSP}
  f_n(\alpha_1,\ldots,\alpha_{n-1}) = t - 2|\alpha| + 1,
\end{equation}
where $|\alpha| = \sum_{i=1}^{n-1} \alpha_i$.
\end{Thm}
\begin{proof}
  Suppose that $A$ has the \SSP. Let
$(\alpha_1,\ldots,\alpha_{n-1}) \in J_{n-1}$. Then $|\alpha| \le
[t/2]$ by Lemma \ref{Lem:SizeOfAlpha}. Hence we have
$(\gin(I):x_n^{t-2|\alpha|})_{|\alpha|} = \gin(I)_{|\alpha|}$. If
$f_n(\alpha_1,\ldots,\alpha_{n-1}) \le t - 2|\alpha|$, then
$x_1^{\alpha_1}\cdots x_{n-1}^{\alpha_{n-1}}x_n^{t-2|\alpha|} \in
\gin(I)$. So $x_1^{\alpha_1}\cdots x_{n-1}^{\alpha_{n-1}} \in
(\gin(I):x_n^{t-2|\alpha|})_{|\alpha|} = \gin(I)_{|\alpha|}$. But
this contradicts that $J_{n-1}$ contains
$(\alpha_1,\ldots,\alpha_{n-1})$. Thus we have
\[
f_n(\alpha_1,\ldots,\alpha_{n-1}) > t - 2|\alpha|. \eqno{(1)}
\]
On the other hand, note that $x_{n-1}^{|\alpha|}x_n^{t-2|\alpha|+1}
\in \gin(I)$ because $R_{t-|\alpha| + 1} = (\gin(I)+
x_n^{t-2|\alpha|+2})_{t-|\alpha| + 1}$. This shows that
$f_n(0,\ldots,0,|\alpha|) \le t-2|\alpha|+1$, by the definition of
$f_n$. Since $\gin(I)$ is strongly stable, we have
\[
f_n(\alpha_1,\ldots,\alpha_{n-1}) \le f_n(0,\ldots,0,|\alpha|) \le
t- 2|\alpha| + 1. \eqno{(2)}
\]
The assertion follows from (1) and (2).

  Conversely, suppose that the equation in (\ref{Eq:SSP}) holds. We need
to show that $(\gin(I):x_n^{t-2i})_i = \gin(I)_i$ and
$x_{n-1}^{i+1}x_n^{t-2i-1} \in \gin(I)_{t-i}$ for all $0 \le i \le
[t/2]$.

  First we will show that $(\gin(I):x_n^{t-2i})_i = \gin(I)_i$ for all $0 \le i \le [t/2]$.
Suppose $x_1^{\alpha_1} \cdots
x_{n-1}^{\alpha_{n-1}}x_n^{i-|\alpha|} \in (\gin(I):x_n^{t-2i})_i$,
where $|\alpha| = \sum_{j=1}^{n-1} \alpha_j \le i$. Then
$x_1^{\alpha_1}\cdots x_{n-1}^{\alpha_{n-1}} x_n^{t-|\alpha|-i} \in
\gin(I)$. Now if $(\alpha_1,\ldots,\alpha_{n-1})$ is an element of
$J_{n-1}$, then $f_n(\alpha_1,\ldots,\alpha_{n-1}) = t - 2|\alpha| +
1$ by the assumption. Hence we have
\[
  t - 2|\alpha| + 1 = f_n(\alpha_1,\ldots,\alpha_{n-1}) \le t
-|\alpha| - i,
\]
by the definition of $f_{n}$. But this contradicts that $|\alpha|
\le i$. So $(\alpha_1,\ldots,\alpha_{n-1}) \notin J_{n-1}$. By
definition of $J_{n-1}$, this implies that either $\alpha_1 \ge f_1$
or there exists $\nu$ such that $1 < \nu \le n-1$,
$(\alpha_1,\ldots,\alpha_{\nu-1}) \in J_{\nu-1}$ and
$f_{\nu}(\alpha_1,\ldots,\alpha_{\nu-1}) \le \alpha_{\nu}$. In any
case $x_1^{\alpha_1} \cdots x_{n-1}^{\alpha_{n-1}} \in \gin(I)$, so
$x_1^{\alpha_1} \cdots x_{n-1}^{\alpha_{n-1}} x_n^{i-|\alpha|} \in
\gin(I)$, and hence $(\gin(I):x_n^{t-2i})_i = \gin(I)_i$.

Finally, to show that $x_{n-1}^{i+1}x_n^{t-2i-1} \in \gin(I)_{t-i}$,
we may assume that $i+1 \le r_1(A)$, otherwise $x_{n-1}^{i+1} \in
\gin(I)$ by the definition of $r_1(A)$, hence the assertion is
achieved. Since $i+1 \le r_1(A) = f_{n-1}(0,\ldots,0)-1$, the
element $(0,\ldots,0,i+1) \in \mathbb{Z}_{\ge 0}^{n-1}$ belongs to
$J_{n-1}$. By the assumption, we have $f_{n}(0,\ldots,0,i+1) =
t-2(i+1)+1 = t-2i-1$, and hence $x_{n-1}^{i+1}x_n^{t-2i-1} \in
\gin(I)_{t-i}$.
\end{proof}

\begin{Coro}\label{Cor:FlagMinGenForSSP}
If $A$ has the \SSP, then
\[
  \{
    x_{n-1}^{r_1(A)+1}, x_{n-1}^{r_1(A)} x_n^{t-2r_1(A)+1}, \ldots,
    x_{n-1}^{i} x_n^{t-2i+1}, \ldots, x_{n-1} x_n^{t-1}, x_n^{t+1}
  \} \subset \G(\gin(I)).
\]
\end{Coro}
\begin{proof}
$x_{n-1}^{r_1(A)+1} \in \G(\gin(I))$ is clear. By Proposition
\ref{Prop:MinimalGen}, $x_{n-1}^{i} x_n^{f_n(0,\ldots,0,i)} \in
\G(\gin(I))$ for all $0 \le i \le r_1(A)$. Since $A$ has the \SSP,
\[
  f_n(0,\ldots,0,i) = t -2i +1,
\]
by Theorem \ref{Thm:SSP}, and hence the assertion follows.
\end{proof}

\begin{Ex}
  Let $I$ and $J$ be the monomial ideal of $R=k[x,y,z,w]$ defined by
\[
I=(x^2, xy, y^3, y^2z, xz^2, yz^2, z^3, \underline{xzw^2 , y^2w^3},
yzw^3, z^2w^3, xw^5, yw^5, zw^5, w^7),
\]
\[
J=(x^2, xy, y^3, y^2z, xz^2, yz^2, z^3, \underline{xzw^3, y^2w^3},
yzw^3, z^2w^3, xw^5, yw^5, zw^5, w^7).
\]
  As shown in Example \ref{Ex:NonUniqueSLP}, the Artinian ring $R/I$
satisfies the \SLP. But since $\max \{i \ | \, (R/I)_i \neq 0 \} =
6$ and $f_4(1,0,1) = 2 < 6 - 2 \times 2 + 1 = 3$, $R/I$ doesn't have
the \SSP. But we can see that $R/J$ satisfies the condition in
Theorem \ref{Thm:SSP} by just reading the $w$-degrees of the minimal
generators of $J$. Hence $R/J$ has the \SSP.
\end{Ex}

\vskip 1cm

\section{\sc Uniquely Determined Minimal Generators of Generic initial ideals}
   We showed that if a standard Artinian $k$-algebra
$A = R/I$ of codimension $n$ has the \WLP, then graded Betti numbers
of $\gin(I)$ are determined just by the Hilbert function of $I$ in
Corollary \ref{Cor:BettiForWLP}. But, if we restrict ourselves to
the case of codimension 3, then we can say much more:
    Throughout this section, we assume that $R=k[x_1,x_2,x_3]$ is a
polynomial ring over a field $k$ and $I$ is a homogeneous Artinian
ideal of $R$. Let $A=R/I=\bigoplus_{i=0}^{t} A_i$, where $t=\max\{i
\ | \, A_i \neq 0 \}$.

    In this section, we show that if $R/I$ has the \WLP, then the minimal generator $T$'s
of $\gin(I)$ with $\max(T) \le 2$ are uniquely determined by the
graded Betti numbers of $R/I$ (Proposition
\ref{Prop:UniquelyDetermineMinimalGeneratorWithWLP}). Together with
Corollary \ref{Cor:BettiForWLP}, this implies that every graded
Betti numbers of $\gin(I)$ can be computed from the graded Betti
numbers of $I$. And if $R/I$ has the \SLP, then we show that the
minimal system of generators of $\gin(I)$ is uniquely determined by
the graded Betti numbers of $R/I$ (Proposition
\ref{Prop:GinSLPUnique}). At last we show that if $R/I$ has the
\SSP, then the minimal system of generators of $\gin(I)$ is uniquely
determined by the Hilbert function of $R/I$ (Proposition
\ref{Prop:GinSSPUnique}). To do so, we use the following theorem
proved by Green.

\begin{Thm}[Cancellation Principle] \cite{Gr}\label{Thm:Cancellation}
  For any homogeneous ideal $J$ in the polynomial ring
$S=k[x_1,\ldots,x_n]$ and any $i$ and $d$, there is a complex of
$k$-modules $V^{d}_{\bullet}$ such that
  \begin{align*}
     V_i^d &\cong \Tor_i^S(\iin(J),k)_d, \\
     H_i(V_{\bullet}^d) &\cong \Tor_i^S(J,k)_d.
  \end{align*}
\end{Thm}

  This theorem means we can obtain the minimal free resolution of
$J$ from the minimal free resolution of $\gin(J)$ by deleting some
adjacent summands of the same degree. Using this theorem, we get the
following lemma.

\begin{Lem}\label{Lem:BettiNumOfGin(I)}
  Suppose that $A=R/I$ has the \WLP.
  If $a = \min\{ d \ | \, \beta_{0,d}(I) \neq 0 \}$, then we have
  \begin{enumerate}
  \item
  $\beta_{0,i}(\gin(I)) =
     \begin{cases} \beta_{1,a+1}(\gin(I)) - \beta_{2,a+2}(\gin(I)) + 1 & \
                       \text{ if } \ i = a,\\
                   \beta_{1,i+1}(\gin(I)) - \beta_{2,i+2}(\gin(I))& \
                       \text{ if } \ i \neq a.
     \end{cases}
  $
  \item For any $i \le r_1(A)+1$,
  \[
  \beta_{1,i+1}(\gin(I)) =
         \beta_{0,i+1}(\gin(I)) + \beta_{1,i+1}(I) -
               \beta_{0,i+1}(I).
  \]
  \end{enumerate}
\end{Lem}
\begin{proof}
  First we will show that the equality in (2) holds. By Theorem
\ref{Thm:Cancellation}, for each $i$, there exist $s_i$, $t_i \in
\mathbb{Z}_{\ge 0}$ such that
\begin{align*}
   \beta_{0,i+1}(\gin(I))& - t_i = \beta_{0,i+1}(I), \\
   \beta_{1,i+1}(\gin(I))& - t_i - s_i = \beta_{1,i+1}(I), \\
   \beta_{2,i+1}(\gin(I))& - s_i = \beta_{2,i+1}(I).
\end{align*}
 By Proposition \ref{Prop:WLP}, if
$i \le r_1(A) + 1$, then $\beta_{2,i+1}(\gin(I)) = 0$, and so $s_i =
0$. This implies that $\beta_{1,i+1}(\gin(I)) - t_i =
\beta_{1,i+1}(I)$, and $\beta_{0,i+1}(\gin(I)) - t_i =
\beta_{0,i+1}(I)$. Hence the equality follows.

   By Eliahou-Kervaire Theorem, we have
   \begin{align*}
       \beta_{1,i+1}(\gin(I))
           = &\sum_{\substack{T \in \G(\gin(I))_{i} \\ \max(T)=2}}
                   \binom{2-1}{1}
              + \sum_{\substack{T \in \G(\gin(I))_{i} \\ \max(T)=3}}
                   \binom{3-1}{1} \\
           = &\sum_{\substack{T \in \G(\gin(I))_{i} \\ \max(T)=1}}
                   1
              + \sum_{\substack{T \in \G(\gin(I))_{i} \\ \max(T)=2}}
                   1
              + \sum_{\substack{T \in \G(\gin(I))_{i} \\ \max(T)=3}}
                   1 \\
              & \ + \sum_{\substack{T \in \G(\gin(I))_{i} \\ \max(T)=3}}
                   1
                  - \sum_{\substack{T \in \G(\gin(I))_{i} \\ \max(T)=1}}
                   1 \\
          =& \begin{cases}
               \beta_{0,a}(\gin(I)) + \beta_{2,a+2}(\gin(I)) - 1 & \ \text{if} \ i = a, \\
               \beta_{0,i}(\gin(I)) + \beta_{2,i+2}(\gin(I)) & \ \text{if} \ i \neq a,
             \end{cases}
   \end{align*}
   and hence the equality in (1) follows.
\end{proof}

\begin{Coro}\label{Cor:UniqueBettiNumofGinForWLP}
  If $A=R/I$ has the \WLP,
then the graded Betti numbers of $\gin(I)$ are completely determined
by the graded Betti numbers of $I$. In particular, if $J$ is a
homogeneous ideal of $R$ such that $R/J$ has the \WLP \ and
$\beta_{i,j}(I) = \beta_{i,j}(J)$ for all $i, j$, then
$\beta_{i,j}(\gin(I)) = \beta_{i,j}(\gin(J))$ for all $i$,$j$.
\end{Coro}
\begin{proof}
  The graded Betti numbers of $I$ determine the Hilbert function of
$A=R/I$. Since $A=R/I$ has the \WLP, $r_1(A)$ is obtained by
Corollary \ref{Cor:r_1WithWLP}. And the graded Betti numbers
$\beta_{2,2+s}(\gin(I))$, for any $s$, and $\beta_{i,i+j}(\gin(I))$,
for $i=0$,$1$ and $j \ge r_1(A) + 2$, are determined by the Hilbert
function of $R/I$, as shown in Corollary \ref{Cor:BettiForWLP}.

Hence it is enough to show that $\beta_{0,j}(\gin(I))$ and
$\beta_{1,j+1}(\gin(I))$ are determined by the graded Betti numbers
of $I$ for $j \le r_1(A) + 1$. Since $\beta_{0,r_1(A)+2}(\gin(I))$
is already known, $\beta_{1,r_1(A)+2}(\gin(I))$ is determined by
Lemma \ref{Lem:BettiNumOfGin(I)} (2), and also
$\beta_{0,r_1(A)+1}(\gin(I))$ is determined by Lemma
\ref{Lem:BettiNumOfGin(I)} (1). Inductively, if we know the value of
$\beta_{0,i+1}(\gin(I))$ for $i \le r_1(A)$, then
$\beta_{1,i+1}(\gin(I))$ and $\beta_{0,i}(\gin(I))$ are determined
by Lemma \ref{Lem:BettiNumOfGin(I)} (2) and (1), respectively. Hence
the assertion follows.
\end{proof}

Now we will show that if $A=R/I$ has the \WLP, then the minimal
generator $T$ of $\gin(I)$ with $\max(T) \le 2$ is determined by the
graded Betti numbers of $I$. We need the following easy lemma.

\begin{Lem}\label{Lem:UniqueSequence}
  Let $a_1,\ldots,a_n$ be a non-decreasing sequence of integers.
Let $\beta_d := \sharp \{ i \ | \, a_i = d \}$ for each $d \in
\mathbb{Z}$. If $\beta_d$ is given for all $d$, then the integers
$a_1,\ldots, a_n$ are determined by $\beta_d$. $\qed$
\end{Lem}

\begin{Prop}\label{Prop:UniquelyDetermineMinimalGeneratorWithWLP}
  Suppose that $I, J$ are homogeneous Artinian ideals of $R$ having the same
graded Betti numbers. If $R/I$ and $R/J$ have the \WLP, then $\{ T
\in \G(\gin(I)) \ | \, \max(T) \le 2 \} = \{ T \in \G(\gin(J)) \ |
\, \max(T) \le 2 \}$. In particular, $\gin(I)_d = \gin(J)_d$ for any
$d \le r_1(A)$.
\end{Prop}
\begin{proof}
  By Corollary \ref{Cor:UniqueBettiNumofGinForWLP}, we have to show
that the minimal generator $T$'s of $\gin(I)$ with $\max(T) \le 2$
are determined by the graded Betti numbers of $\gin(I)$. At first
note that the minimal generator $T$'s of $\gin(I)$ such that
$\max(T) \le 2$ are of the forms;
\[
 x_1^a, x_1^{a-1}x_2^{f_2(a-1)}, x_1^{a-2}x_2^{f_2(a-2)}, \ldots, x_1x_2^{f_2(1)},
x_2^{f_2(0)},
\]
where $a := \min \{i \ | \, \beta_{0,i}(I) \neq 0 \}$. So it is
enough to show that each $f_2(i)$ is determined by the graded Betti
numbers of $\gin(I)$. Since $\gin(I)$ is strongly stable,
\[
 \begin{array}{rl}
  a \ \le \  a-1 + f_2(a-1) &\le \ a - 2 + f_2(a-2) \ \le \  \cdots \\
    &\le \ 1 + f_2(1) \ \le \ f_2(0) = r_1(A)+1.
 \end{array}
\]
And note that for $i \le r_1(A)$, $\beta_{0,i}(\gin(I))$ is the
number of minimal generator $T$'s of $\gin(I)$ such that $\max(T)
\le 2$ and $\deg(T) = i$, by Proposition \ref{Prop:WLP}. This
implies that $\beta_{0,i}(\gin(I)) = \sharp \{ j \ | \, j + f_2(j) =
i \}$. Using Lemma \ref{Lem:UniqueSequence}, we can see that $f_2$
is determined by $\beta_{0,i}(\gin(I))$.

And the last assertion follows by Proposition \ref{Prop:WLP}.
\end{proof}

  If $R/I$ has the \SLP, then we will show that the minimal
system of generators of $\gin(I)$ is determined by the graded Betti
numbers of $I$. Remind that determining the minimal system of
generators of $\gin(I)$ is equivalent to computing $f_1$, $f_2$, and
$f_3$ defined as in (\ref{f's}). From Proposition
\ref{Prop:UniquelyDetermineMinimalGeneratorWithWLP}, we already know
that $f_1$, $f_2$ are determined by the graded Betti numbers of $I$.
Hence it is enough to show that $f_3$ is also determined by those of
$I$ in the case that $A=R/I$ has the \SLP. Now
\begin{equation} \label{Set:IndexMinGenWithn=3}
\begin{split}
  J_2 = \left\{(\alpha_1,\alpha_{2})
              \left| \begin{array}{l}
                       0 \le \alpha_1 < f_1 = a,  \\
                       0 \le \alpha_2 < f_2(\alpha_1)
                     \end{array}
              \right.
        \right\},
\end{split}
\end{equation}
where $f_1=a$ and $f_2$ is as shown in the proof of Proposition
\ref{Prop:UniquelyDetermineMinimalGeneratorWithWLP}. Partition
$J_{2}$ into the sets $J[t]$
           defined by
 \[
    J[t] := \{ (\alpha_1, \alpha_{2}) \in J_{2} \ | \,
                \alpha_1 + \alpha_2 = t \},
 \]
for $0 \le t \le r_1(A)$. Note that $J_{2} = \cup_{t=0}^{r_1(A)}
J[t]$, and the minimal generator of $\gin(I)$ corresponding to the
element $(\alpha_1, \alpha_2) \in J[t]$ is
$x_1^{\alpha_1}x_2^{\alpha_2}x_3^{f_3(\alpha_1,\alpha_2)}$, whose
degree is $\alpha_1 + \alpha_2 + f_3(\alpha_1,\alpha_2) =
t+f_3(\alpha_1,\alpha_2)$. We will say that $\alpha_1 + \alpha_2 +
f_3(\alpha_1,\alpha_2)$ is the degree of $(\alpha_1,\alpha_2)$.

\begin{Prop}\label{Prop:GinSLPUnique}
  If $R/I$ has the \SLP, then the minimal system of
generators of $\gin(I)$ is determined by the graded Betti numbers of
$I$. In particular, if $J$ is a homogeneous ideal of $R$ such that
$\beta_{i,j}(J) = \beta_{i,j}(I)$ for all $i, j$, and if $R/J$ has
the \SLP, then $\gin(I)=\gin(J)$.
\end{Prop}
\begin{proof}
  Let $(\alpha_1,\alpha_2)$, $(\beta_1,\beta_2) \in J[t]$. If
$\alpha_1 \le \beta_1$, then
\[
   f_3(\alpha_1,\alpha_2) \ge f_3(\beta_1,\beta_2),
\]
by Remark \ref{Cond:ElementInJ_i} (3), and hence
\begin{align*}
  \alpha_1 + \alpha_2 + f_3(\alpha_1,\alpha_2) & = t + f_3(\alpha_1,\alpha_2) \\
              & \ge t + f_3(\beta_1,\beta_2) \\
              &= \beta_1 + \beta_2 + f_3(\beta_1,\beta_2).
\end{align*}
This shows that, for each $t$, we have a nondecreasing sequence
$\{a_{t,i}\}$ of degrees of elements in $J[t]$: %
Indeed, suppose that $J[t] = \{ (\alpha_{1,1},\alpha_{1,2}), \ldots,
(\alpha_{j_t,1},\alpha_{j_t,2}) \}$ and $\alpha_{j,1} \ge
\alpha_{j+1,1}$. If we set $a_{t,i} = t +
f_3(\alpha_{i,1},\alpha_{i,2})$, then we have
\[
   a_{t,1} \le a_{t,2} \le \cdots \le a_{t,|J[t]|},
\]
where $|J[t]|$ is the number of elements in $J[t]$. Note that
$a_{t,|J[t]|} = t + f_3(0,t)$, as shown in the above. And since
$A=R/I$ has the \SLP, $(\alpha_1,\alpha_2) \in J[t]$ implies that
\[
 f_{3}(0,t+1) + 1 \le f_{3}(\alpha_1,\alpha_{2}),
\]
 and hence
\[
t + 1 + f_{3}(0,t+1) \le t + f_{3}(\alpha_1,\alpha_{2}),
\]
by Theorem \ref{Thm:SLP}. This shows that $a_{t+1,|J[t+1]|} \le
a_{t,1}$ for any $0 \le t \le r_1(A) - 1$. Thus we have a
non-decreasing sequence of degrees of the elements in the whole set
$J_2$:
\[
  a_{r_1(A),1} \le a_{r_1(A),2} \le \ldots \le a_{r_1(A),|J[r_1(A)]|}
   \le a_{r_1(A)-1,1} \le \ldots \le a_{1,|J[1]|} \le a_{0,1}.
\]
Since $\beta_{2,2+d}(\gin(I))$ is the number of the minimal
generator $T$'s of $\gin(I)$ such that $\deg(T) = d$ and $\max(T) =
3$, i.e.,
\begin{align*}
    \beta_{2,2+d}(\gin(I)) &= \sharp \{ (\alpha_1,\alpha_2) \in J_2 \ |
\, \alpha_1 + \alpha_2 + f_3(\alpha_1,\alpha_2) = d \} \\
                           &= \sharp \{ (i,j) \ | \, a_{i,j} = d \},
\end{align*}
${a_{i,j}}$ is uniquely determined by $\beta_{2,2+d}(\gin(I))$, by
Lemma \ref{Lem:UniqueSequence}, so we are done.
\end{proof}

  Until now, we show that if $A=R/I$ is a standard Artinian graded
$k$-algebra of codimension 3 which has the \SLP, then the minimal
system of generators of $\gin(I)$ is determined by the graded Betti
numbers of $I$. But if $A$ has the \SSP, then we will show that the
minimal system of generators of $\gin(I)$ is determined by the
Hilbert function of $A$.

 For any $(\alpha_1,\alpha_2)$,
$(\beta_1,\beta_2) \in \mathbb{Z}_{\ge 0}^{2}$, we say that
$(\beta_1,\beta_2) \le (\alpha_1,\alpha_2)$ if $x_1^{\beta_1}
x_2^{\beta_2} \le x_1^{\alpha_1} x_2^{\alpha_2}$ with respect to the
reverse lexicographic order. We will use this order in the following
lemma and proposition. And we say that a monomial $T_1 =
x_1^{\alpha_1} x_2^{\alpha_2}$ is obtained from a monomial
$T_2=x_1^{\beta_1} x_2^{\beta_2}$ by elementary Borel move if %
$T_1 = x_1 \frac{T_2}{x_2}$. In that case, we also say that
$(\alpha_1,\alpha_2)$ is obtained from $(\beta_1,\beta_2)$ by
elementary Borel move.

 Set $|\alpha| := \alpha_1 + \alpha_2$. Note that if $|\alpha| = |\beta|$
and $(\beta_1,\beta_2) \le (\alpha_1,\alpha_2)$, then
$x_1^{\alpha_1} x_{2}^{\alpha_{2}}$ is obtained from $x_1^{\beta_1}
x_{2}^{\beta_{2}}$ by consecutive elementary Borel move.

\begin{Lem}\label{Lem:Contin_of_minimal_Generator}
Suppose that $x_1^{\alpha_1} x_{2}^{\alpha_{2}}x_{3}^{\gamma}$,
$x_{2}^{|\alpha|}x_3^{\gamma}$ are minimal generators of $\gin(I)$
for some $(\alpha_1,\alpha_{2}) \in \mathbb{Z}_{\ge 0}^{2}$. Then
$x_1^{\beta_1} x_{2}^{\beta_{2}} x_{3}^{\gamma}$ is also a minimal
generators of $\gin(I)$, for any $(\beta_1,\beta_{2})$ $\le$
$(\alpha_1,\alpha_{2})$ with $|\beta| = |\alpha|$.
\end{Lem}
\begin{proof}
Set $T_{\alpha}:=x_1^{\alpha_1} x_{2}^{\alpha_{2}}x_{3}^{\gamma}$,
and $T_{\beta}:=x_1^{\beta_1} x_{2}^{\beta_{2}}x_{3}^{\gamma}$. By
induction, it is enough to show that if $(\alpha_1,\alpha_2)$ is
obtained from $(\beta_1,\beta_2)$ by elementary Borel move, then
$T_{\beta} \in \G(\gin(I))$. Note that, in that case,
$(\alpha_1,\alpha_2) = (\beta_1+1,\beta_2-1)$. Since
$x_{2}^{|\beta|}x_{3}^{\gamma} \in \G(\gin(I))$ and $\gin(I)$ is
strongly stable, $T_{\beta} \in \gin(I)$.  So it is enough to show
that $T_{\beta}$ is a minimal generator. If not, there exists
$(\delta_1,\delta_{2},\delta) \in \mathbb{Z}_{\ge 0}^{3}$ such that
$\delta \le \gamma$, $\delta_l \le \beta_l$ for any $l = 1, 2$, and
$T_{\delta} := x_1^{\delta_1} x_{2}^{\delta_{2}} x_{3}^{\delta} \in
\gin(I)$. If $\delta_2 < \beta_2$, then $\alpha_2 = \beta_2 - 1 \ge
\delta_2$, hence $T_{\delta}$ divides $T_{\alpha}$, this contradicts
that $T_{\alpha}$ is a minimal generator of $\gin(I)$. Thus we have
$\delta_2 = \beta_2$. But, in this case, $T'_{\delta} = x_1
\frac{T_{\delta}}{x_2} \in \gin(I)$ divides $T_{\alpha}$, a
contradiction. This shows that $T_{\beta} \in \G(\gin(I))$.
\end{proof}

\begin{Prop}\label{Prop:GinSSPUnique}
  If $A=R/I$ has the \SSP, then the minimal system generators of $\gin(I)$
is determined by the Hilbert function of $I$. In particular, if $J$
is a homogeneous ideal of $R$ with $H(R/I,d)=H(R/J,d)$ for all $d$,
and if $R/J$ has the \SSP, then $\gin(I) = \gin(J)$.
\end{Prop}
\begin{proof}
  Let $t = \max \{ i \ | \, A_i \neq 0 \}$. Since $A$ has the
\SSP, $A$ has also the \WLP. From Corollary \ref{Cor:r_1WithWLP} and
Corollary \ref{Cor:BettiForWLP}, $r_1(A)$ and
$\beta_{2,2+d}(\gin(I))$ are given by the Hilbert function of $R/I$.

For $0 \le d \le r_1(A)$, let $B_d := \{ (\alpha_1,\alpha_{2}) \in
\mathbb{Z}_{\ge 0}^{2} \ | \, \alpha_1 + \alpha_2 = d \}$, and let
$J[d]$ be the set consisting of the smallest
$\beta_{2,2+(t-d+1)}(\gin(I))$ elements in the set $B_d$ with
respect to the order we defined. We will show that $J_{2}$ is the
disjoint union of $J[d]$, $0 \le d \le r_1(A)$. Since disjointness
is clear by the definition of $J[d]$, it is enough to show that
$J_2$ is the union of the sets $J[d]$.

Let $(\alpha_1,\alpha_{2}) \in J_{2}$. Then $x_1^{\alpha_1}
x_{2}^{\alpha_{2}} x_3^{f_3(\alpha_1,\alpha_{2})} \in \G(\gin(I))$,
by Remark \ref{Cond:ElementInJ_i} (1). Since $A$ has the \SSP,
\[
  f_3(\alpha_1,\alpha_{2}) = t-2|\alpha| + 1,
\]
where $|\alpha| = \alpha_1 + \alpha_2$, by Theorem \ref{Thm:SSP}.
Note that the degree of $x_1^{\alpha_1} x_{2}^{\alpha_{2}}
x_3^{t-2|\alpha| + 1}$ is $t - |\alpha| + 1$ and $|\alpha| \le
r_1(A)$. Otherwise, $x_1^{\alpha_1}
x_{2}^{\alpha_{2}}x_3^{t-2|\alpha| + 1}$ cannot be a minimal
generator of $\gin(I)$, since $x_{2}^{r_1(A)+1} \in \gin(I)$ and
$\gin(I)$ is strongly stable. Hence $x_{2}^{|\alpha|}
x_3^{t-2|\alpha|+1} \in \G(\gin(I))$, by Corollary
\ref{Cor:FlagMinGenForSSP}. Suppose that $(\alpha_1,\alpha_{2})$ is
the $l$-th smallest element of $B_{|\alpha|}$. By Lemma
\ref{Lem:Contin_of_minimal_Generator}, there are at least $l$
minimal generator $T$'s of $\gin(I)$ such that $\max(T) = 3$ and
$\deg(T) = t-|\alpha| + 1$. So $l \le
\beta_{2,2+(t-|\alpha|+1)}(\gin(I))$, and hence
$(\alpha_1,\alpha_{2}) \in J[|\alpha|]$. This shows that
$\sum_{d=0}^{r_1(A)} \sharp J[d] \ \ge \sharp J_{2}$. And
\begin{align*}
  \sharp J_{2} & = \sharp \{ T \in \G(\gin(I) \ | \, \max(T) = 3 \} \\
                 & = \sum_{i=0}^{t+1} \beta_{2,2+i}(\gin(I)) \
                   \ge \sum_{d=0}^{r_1(A)} \sharp J[d] \ \ge \sharp
J_{2},
\end{align*}
by Remark \ref{Cond:ElementInJ_i} (1). Hence $J_2$ is the disjoint
union of the sets $J[d]$, $0 \le d \le r_1(A)$. The assertion
follows by Remark \ref{Remk:JDetermineGin}.
\end{proof}

\begin{Coro}\label{Cor:NSCondForSSP}
 Suppose that $A=R/I$ has the \SSP. If $J$ is a homogeneous ideal of $R$
having the same Hilbert function with $I$, and if $R/J$ has the
\WLP, then $R/J$ has the \SSP \ if and only if
\[
   \{ x_{2}^{r_1(A)+1}, x_{2}^{r_1(A)} x_3^{t-2r_1(A)+1},
\ldots, x_{2}^{i}x_3^{t-2i+1}, \ldots, x_{2} x_3^{t-1}, x_3^{t+1} \}
\subset \gin(J),
\]
where $t = \max \{ i \ | \, A_i \neq 0 \}$.
\end{Coro}
\begin{proof}
$(\Rightarrow)$ It follows from Proposition \ref{Prop:GinSSPUnique}
and Corollary \ref{Cor:FlagMinGenForSSP}.

$(\Leftarrow)$ Note that every element of the left-hand side is a
minimal generator of $\gin(I)$, as shown in the proof of Proposition
\ref{Prop:GinSSPUnique}. Because the other minimal generator $T$'s
of $\gin(I)$ with $\max(T)=3$ are obtained from $x_{2}^{i}
x_3^{t-2i+1}$ by consecutive elementary Borel move for some $0 \le i
\le r_1(A)$, we have $\{ T \in \G(\gin(I)) \ | \, \max(T) = 3 \}
\subset \gin(J)$. Since $H(R/I,d)=H(R/J,d)$, $\dim_k \gin(I)_d =
\dim_k \gin(J)_d$ for all $d$. This shows that $\{ T \in \G(\gin(I))
\ | \, \max(T) = 3 \} = \{ T \in \G(\gin(J)) \ | \, \max(T) = 3 \}$,
by Proposition \ref{Prop:WLP}. Hence the assertion follows from
Remark \ref{Remk:JDetermineGin}.
\end{proof}

  Proposition \ref{Prop:GinSLPUnique} and Proposition
\ref{Prop:GinSSPUnique} do not work if $A=R/I$ is a standard graded
Artinian $k$-algebra of codimension $\ge 4$, as shown in the
following example.

\begin{Ex}\label{Ex:NonUniqueSSP}
Consider two strongly stable ideals of $R=k[x,y,z,w]$ defined as
\begin{align*}
I=(&x^2, xy^2, y^4, \underline{y^3z, xyz^3, y^2z^3, xz^5}, yz^5, z^7, \\
   &z^6w, \underline{xz^4w^3}, yz^4w^3, z^5w^3, \underline{xyz^2w^5},
        y^2z^2w^5, xz^3w^5, yz^3w^5, z^4w^5, \\
   &y^3w^7, xyzw^7, y^2zw^7, xz^2w^7, yz^2w^7, z^3w^7,
        xyw^9, y^2w^9, xzw^9, yzw^9, z^2w^9, \\
   &xw^{11}, yw^{11}, zw^{11}, w^{13}), \\
I'=(&x^2, xy^2, y^4, \underline{xyz^2, y^3z^2, xz^4, y^2z^4}, yz^5, z^7, \\
   &z^6w, \underline{y^2z^3w^3}, yz^4w^3, z^5w^3, \underline{y^3zw^5},
        y^2z^2w^5, xz^3w^5, yz^3w^5, z^4w^5, \\
   &y^3w^7, xyzw^7, y^2zw^7, xz^2w^7, yz^2w^7, z^3w^7,
        xyw^9, y^2w^9, xzw^9, yzw^9, z^2w^9, \\
   &xw^{11}, yw^{11}, zw^{11}, w^{13}).
\end{align*}
  Note that $I$ and $I'$ have the same graded Betti numbers. And for
any $(\alpha_1,\alpha_2,\alpha_3) \in J_3$,
$f_4(\alpha_1,\alpha_2,\alpha_3) = t - 2|\alpha| + 1$, where
$t=\max\{ i \ | \, (R/I)_i \neq 0 \} = \max\{ i \ | \, (R/I')_i \neq
0 \} = 12$. Hence they have the \SSP.
\end{Ex}\vskip 1cm

\bibliographystyle{amsalpha}

\end{document}